\documentclass{lms}
\usepackage{amsmath,amsfonts,amssymb,url}

\newcommand{\st}{s.t.}

\DeclareFontFamily{OT1}{rsfs}{}
\DeclareFontShape{OT1}{rsfs}{n}{it}{<-> rsfs10}{}
\DeclareMathAlphabet{\mathscr}{OT1}{rsfs}{n}{it}

\DeclareMathOperator{\num}{num}
\DeclareMathOperator{\den}{den}
\DeclareMathOperator{\lcm}{lcm}

\DeclareMathOperator{\rad}{rad}
\DeclareMathOperator{\mo}{\,mod}
\DeclareMathOperator{\sq}{sq}

\DeclareMathOperator{\Area}{Area}

\newtheorem{prop}{Proposition}[section]
\newtheorem{thm}[prop]{Theorem}
\newtheorem{cor}[prop]{Corollary}
\newtheorem{lem}[prop]{Lemma}
\newnumbered{defn}{Definition}
\newenvironment{Rem}{{\bf Remark.}}{}
\numberwithin{equation}{section}
\title{The parity problem for reducible cubic forms}
\author{H. A. Helfgott}
\classno{11N32 (primary), 11N35, 11N36 (secondary)}
\begin{document}
\maketitle
\begin{abstract}
Let $f\in \mathbb{Z}\lbrack x,y\rbrack$ be a reducible homogeneous
polynomial of degree $3$. We show that $f(x,y)$ has an even
number of prime factors as often as an odd number of prime factors.
\end{abstract}
\section{Introduction}
The Liouville function $\lambda(n)$ is defined on
the set of non-zero rational integers as follows:
\begin{equation}\label{eq:deflam}
\lambda(n) = \prod_{p|n} (-1)^{v_p(n)} .
\end{equation}
We will find it convenient to choose a value for $\lambda(0)$; we adopt the
convention that $\lambda(0) = 0$. 

Let $f\in \mathbb{Z}\lbrack x,y\rbrack$ be a homogeneous polynomial not
of the form $c\cdot g^2$, $c\in \mathbb{Z}$, $g\in \mathbb{Z}\lbrack x,y\rbrack$. Then, it is believed,
\begin{equation}\label{eq:honeysuckle}
\lim_{N\to \infty} \frac{1}{N^2} \sum_{-N\leq x,y\leq N} \lambda(f(x,y)) = 0 .
\end{equation}
This conjecture can be traced to Chowla (\cite{Ch}, p. 96); it is closely 
related to the Bunyakovsky/Schinzel conjecture on primes represented by
irreducible polynomials. 

The one-variable analogue of (\ref{eq:honeysuckle}) is classical 
for $\deg f = 1$ and quite hopeless for $\deg f> 1$. We know 
(\ref{eq:honeysuckle}) itself when $\deg f\leq 2$.
(The main ideas of the proof go back to de la Vall\'ee-Poussin
(\cite{DVP1}, \cite{DVP2}); see \cite{He2}, \S 3.3, for an exposition.)
The problem of proving (\ref{eq:honeysuckle}) when $\deg f\geq 3$ has
remained open until now: 
sieving is forestalled by the parity problem (\cite{Se}), which 
Chowla's conjecture may be said to embody in its pure form.

We prove (\ref{eq:honeysuckle}) for $f$ reducible of degree $3$.
In a companion paper (\cite{Heirr}), we prove (\ref{eq:honeysuckle})
for $f$ irreducible of degree $3$.

%We will show that the average of
%$\lambda(f(x,y))$ is zero not just over $\mathbb{Z}^2$, but, more
%generally, over the intersection of any lattice coset and any convex set.

Part of the importance of Chowla's conjecture resides in its applications 
to problems of parity outside analytic number theory.
Knowing that (\ref{eq:honeysuckle}) holds for $\deg f = 3$ allows us to
conclude that in certain one-parameter families of elliptic curves the
root number $W(E) = \pm 1$ averages to 
$0$ (\cite{He}, Proposition 5.6). In \S \ref{sec:theman}, we will show that the 
two-parameter family $y^2 = x (x+a) (x + b)$ has average root number $0$
as well. In the process, we will see that, for some $f$,
(\ref{eq:honeysuckle}) is robust under certain twists by characters
to variable moduli.

\section{Preliminaries}
\subsection{Anti-sieving}
In the next two lemmas we use an upper-bound sieve not to find almost-primes, but to split the integers multiplicatively, with the almost-primes as an error term. A treatment
by means of a cognate of Vaughan's identity would also be possible, but much more cumbersome.
The error term would be the same.
%Please check the constants in this lemma!
\begin{lem}\label{lem:sieve}
Let $\mathscr{P}\subset \{M_1,M_1+1,\dotsc,M_2-1\}$ 
be a set of primes between the positive integers
$M_1$ and $M_2$.
Then there are $\sigma_d\in \mathbb{R}$ with $|\sigma_d|\leq 1$
and support on 
\[\{M_1\leq d < M_2 : p|d \Rightarrow p\in \mathscr{P}\}\]
such that for any $a$, $m$, $N_1$ and $N_2$ with $0< m<M_1$ and
$N_2\geq N_1$,
\[\mathop{\sum_{N_1\leq n< N_2}}_{n\equiv a \mo m} 
\left|1-\sum_{d|n} \sigma_d\right| \ll 
\mathop{\prod_{M_1\leq p< M_2}}_{p\notin \mathscr{P}} 
\left(1 - \frac{1}{p}\right)^{-1}
\cdot
\frac{\log M_1}{\log M_2}
\frac{N_2-N_1}{m}
 + M_2,\]
where the implied constant is absolute.
\end{lem}
\begin{proof}
We shall follow the nomenclature in \cite{IK}, p. 159. We let 
\[\kappa = 1,\;\;\;\;\; y = M_2,\;\;\;\;\; g(d) = \begin{cases} 1/d &\text{if $p|d \Rightarrow
p\in \mathscr{P}$,}\\0 &\text{otherwise.}\end{cases}\] 
Then, by Fundamental Lemma 6.3 in \cite{IK}, there is a sequence
of real numbers $(\lambda_d^+)$ such that
\[\lambda_d^+ = 1,\;\;\;\; |\lambda_d^+|\leq 1 \text{ for all $d$},\;\;\;\;
\lambda_d^+ = 0 \text{ if $d\geq y$},\;\;\;\;\sum_{d|n} \lambda_d^+ \geq 0 \text{ for every $n$,}\]
and
\begin{equation}\label{eq:cure}
\sum_{d|P(z)} \lambda_d^+ g(d) \ll \prod_{p<z} (1 - g(p)),\end{equation}
where $z = M_2$ and $P(z) = \prod_{p<z} p$. Note now that 
\[\prod_{p<z} (1 - g(p)) \ll \frac{\log M_1}{\log M_2}
\prod_{M_1\leq p<M_2 : p\notin \mathscr{P}} (1 - 1/p)^{-1}.\] Hence
\[\begin{aligned}
%\mathop{\sum_{N_1\leq n < N_2}}_{n\equiv a \mo m} \left|
%\mathop{\sum_{d|n}}_{p|d \Rightarrow p\in \mathscr{P}} \lambda_d^+
%\right| &=
 \mathop{\sum_{N_1\leq n < N_2}}_{n\equiv a \mo m}
\mathop{\sum_{d|n}}_{p|d \Rightarrow p\in \mathscr{P}} \lambda_d^+
&= \mathop{\sum_d}_{p|d \Rightarrow p\in \mathscr{P}} \lambda_d^+
 \cdot \left(\frac{N_2 - N_1}{m d} + O(1)\right)\\
&= 
\frac{N_2 - N_1}{m} \sum_{d|P(z)} \lambda_d^+ g(d) \;\;+
\mathop{\sum_d}_{p|d \Rightarrow p\in \mathscr{P}} O(\lambda_d^+) 
\\
&\ll \mathop{\prod_{M_1\leq p < M_2}}_{p\notin \mathscr{P}}
\left( 1 - \frac{1}{p}\right)^{-1} \cdot \frac{\log M_1}{\log M_2}
\frac{N_2 - N_1}{m} + M_2,\end{aligned}\]
where all implied constants are absolute.
We set \[\sigma_1 = 0,\;\;\;\;\;\;\;\;\;\;\;\;\;\sigma_d =
\begin{cases} -\lambda_d^+ &\text{if $p|d \Rightarrow p\in \mathscr{P}$}\\
0 &\text{otherwise}\end{cases}\; \text{ for $d\neq 1$.}\]
Since $\sum_{d|n} \lambda_d^+\geq 0$ for every $n$, we are done.
\end{proof}
\begin{Rem}
Fundamental Lemma 6.3 in \cite{IK} employs the Rosser-Iwaniec
sieve, and thus gives an optimized bound for the constant in (\ref{eq:cure}).
As any constant would do for our purposes, we could use somewhat weaker
 results, such as Brun's 1920 sieve (\cite{Gr}, \S 3.4). All the same, we are 
using -- and need -- a result different from some that go by the name 
of ``fundamental lemma''
in the older literature: we are not assuming that $\log M_2 = o(\log N)$,
and we are not requiring asymptotics.
\end{Rem}
\begin{lem}\label{sieve2}
Let $K/\mathbb{Q}$ be a number field.
Let $\jmath:K\to \mathbb{R}^{\deg (K/\mathbb{Q})}$ be a bijective
$\mathbb{Q}$-linear map
taking $\mathscr{O}_K$ to $\mathbb{Z}^{\deg (K/\mathbb{Q})}$.
Let $\mathscr{P}$ be a set of prime ideals of $K$
whose norms lie between the positive integers
$M_1$ and $M_2$.
Then there are $\sigma_{\mathfrak{d}}
\in \mathbb{R}$ with $|\sigma_{\mathfrak{d}}|\leq 1$
and support on
\begin{equation}\label{eq:schi}
\{\mathfrak{d} : M_1\leq N \mathfrak{d} < M_2,\;
\mathfrak{p} | \mathfrak{d} \Rightarrow \mathfrak{p}\in \mathscr{P}\}
\end{equation}
such that for any positive integer $N$, any lattice
coset $L\subset \mathbb{Z}^{\deg(K/\mathbb{Q})}$ with index 
$\lbrack \mathbb{Z}^{\deg(K/\mathbb{Q})} : L\rbrack<M_1$ and
any convex set 
  $S\subset \lbrack -N,N\rbrack^{\deg(K/\mathbb{Q})}$,
\[\sum_{\jmath(x)\in S\cap L} 
\left|1-\sum_{\mathfrak{d}|x} \sigma_{\mathfrak{d}} \right| \ll
\mathop{\prod_{M_1\leq p< M_2}}_{p\notin \mathscr{P}} \left(1 - \frac{1}{p}\right)^{-1}
\cdot
\frac{\log M_1}{\log M_2}
\frac{\Area(S)}{\lbrack \mathscr{O}_K:L\rbrack}
 + N^{\deg(K/\mathbb{Q})-1} M_2^2,\]
where 
the implied constant depends only on $K$.
\end{lem}
\begin{proof}
Set $\lambda_{\mathfrak{d}}^+$ as in a generalized
 Rosser--Iwaniec sieve 
(\cite{Col2}) with sieving set $\mathscr{P}$
and upper cut $z=M_2$. Proceed as in the proof of Lemma 
\ref{lem:sieve}.
Set $\sigma_{(1)} = 0$, $\sigma_\mathfrak{d}
= -\lambda_\mathfrak{d}^+$ for $\mathfrak{d}\ne (1)$.
\end{proof}

\subsection{Extensions of the Liouville function}
We define $\lambda$ on $\mathbb{Q}$ by
\begin{equation}\label{eq:extrat}
\lambda\left(\frac{n_0}{n_1}\right) = \frac{\lambda(n_0)}{\lambda(n_1)}
\end{equation}
and on
ideals in a Galois extension $K/\mathbb{Q}$ of degree $n$ by
\begin{equation}\label{eq:exide}
\lambda(\mathfrak{p}_1^{e_1} \mathfrak{p}_2^{e_2} \dotsb \mathfrak{p}_k^{e_k}) =
\prod_{i} \omega^{f(\mathfrak{p_i})\cdot e_i},\end{equation}
where $\omega$ is a fixed $(2n)$th root of unity and
$f(\mathfrak{p_i})$ is the degree of inertia of $\mathfrak{p}_i$ over 
$\mathfrak{p}_i \cap \mathbb{Q}$. Notice that (\ref{eq:exide}) restricts
to (\ref{eq:extrat}), which, in turn, restricts to (\ref{eq:deflam}).
Notice also that the above {\em extension} is different from the natural
{\em generalization} $\lambda_K$:
\begin{equation}\label{eq:exgen}
\lambda_K(\mathfrak{p}_1^{e_1} \mathfrak{p}_2^{e_2} \dotsb \mathfrak{p}_k^{e_k}) =
\prod_{i} (-1)^{e_i}.\end{equation}

\subsection{Quadratic forms}
We will consider only quadratic forms $a x^2 + b x y + c y^2$ with 
integer coefficients $a,b,c\in \mathbb{Z}$. A quadratic form $a x^2 + b x y
+ c y^2$ is {\em primitive} if $\gcd(a,b,c)=1$.  

Let $n$ be a rational integer. We denote by $\sq(n)$ the largest 
positive integer whose square divides $n$. Define
\[d_n = \begin{cases} \sq(n) &\text{if $4\nmid n$}\\
\sq(n)/2 &\text{if $4|n$.}\end{cases}\]

\begin{lem}\label{lem:quad}
Let $Q(x,y) = a x^2 + b x y + c y^2$ be a primitive, irreducible quadratic form. Let $K = \mathbb{Q}(\sqrt{b^2 - 4 a c})$. Then there are 
algebraic integers $\alpha_1, \alpha_2\in \mathscr{O}_K$ linearly
independent over $\mathbb{Q}$
such that 
\[ Q(x,y) = \frac{N(x \alpha_1 + y \alpha_2)}{a}\]
for all $x,y\in \mathbb{Z}$. The subgroup $\mathbb{Z} \alpha_1 + \mathbb{Z} \alpha_2$
of $\mathscr{O}_K$ has index $\lbrack \mathscr{O}_K : \mathbb{Z} \alpha_1 +
\mathbb{Z} \alpha_2\rbrack = d_{b^2 - 4 a c}$.
\end{lem}
\begin{proof}
Set $\alpha_1 = a$, $\alpha_2 = \frac{b + \sqrt{b^2 - 4 a c}}{2}$. 
\end{proof}

\subsection{Lattices and convex sets}

A \emph{lattice} is a subgroup of $\mathbb{Z}^n$ of finite index; a \emph{%
lattice coset} 
 is a coset of such a subgroup. By the {\em index} of a lattice
coset we mean the index of the lattice of which it is a coset.
For any lattice cosets $L_1$, $%
L_2$ with $\gcd(\lbrack \mathbb{Z}^n : L_1 \rbrack , \lbrack \mathbb{Z}^n :
L_2 \rbrack) = 1$, the intersection $L_1\cap L_2$ is a lattice coset with 
\begin{equation}
\lbrack \mathbb{Z}^n : L_1\cap L_2\rbrack = \lbrack \mathbb{Z}^n :
L_1\rbrack \lbrack \mathbb{Z}^n : L_2\rbrack .
\end{equation}
In general, if $L_1$, $L_2$ are lattice cosets, then $L_1\cap L_2$ is either
the empty set or a lattice coset such that 
\begin{equation}  \label{eq:intersl} 
\lcm(\lbrack \mathbb{Z}^n : L_1\rbrack,\lbrack \mathbb{Z}^n : L_2\rbrack) 
\mid
\lbrack \mathbb{Z}^n : L_1\cap L_2\rbrack,\;\;\; \lbrack \mathbb{Z}^n :
L_1\cap L_2\rbrack \mid \lbrack \mathbb{Z}^n : L_1\rbrack \lbrack \mathbb{Z}^n
: L_2\rbrack .
 \end{equation}
Since $\mathbb{Z}^n/L_j$ ($j=1,2$) is a quotient of $\mathbb{Z}^n/(L_1\cap L_2)$,
we must have $\lbrack \mathbb{Z}^n : L_j\rbrack | \lbrack \mathbb{Z}^n :
L_1 \cap L_2\rbrack$. The first property in (\ref{eq:intersl}) follows. 
Two distinct elements
of $\mathbb{Z}^n / (L_1 \cap L_2)$ cannot be congruent modulo both $L_1$ and
$L_2$. Thus, the natural map $\mathbb{Z}^n/(L_1 \cap L_2) \mapsto
\mathbb{Z}^n/L_1 \times \mathbb{Z}^n/L_2$ must be injective. The second
property in (\ref{eq:intersl}) follows.

For $S\subset \lbrack -N,N\rbrack^n$ a convex set and $L\subset \mathbb{Z}^n$ a lattice
coset, 
\begin{equation}  \label{eq:wbrl}
\#(S\cap L) = \frac{\Area(S)}{\lbrack \mathbb{Z}^n :L\rbrack} + O(N^{n-1}),
\end{equation}
where the implied constant depends only on $n$. One can prove (\ref{eq:wbrl})
easily: slice $S$ and $L$ by hyperplanes and use induction on $n$.

\subsection{Linear bounds}
Landau showed (\cite{La}) that
there is an effective constant $c>0$ such that, for every $k\geq 0$, there
is at most one primitive character $\chi$
of modulus $q\in \{2^{2^k},2^{2^k}+1,\dotsc,
2^{2^{k+1}}-1\}$ such that $L(s,\chi)$ has an exceptional (``Siegel'') zero 
$\beta>1 - c/\log q$. We call such a modulus $q$ {\em exceptional}.
 By Siegel's methods 
(\cite{Si}; vd.\ also \cite{Pr}, p. 74--75), it follows that
\begin{equation}\label{eq:sw}
\left|\mathop{\sum_{n\leq x}}_{n\equiv a \mo m} \lambda(n)\right| \ll
x e^{-C \sqrt{\log x}}\end{equation}
for any $x\geq 1$ and any $a$, $m$ such that $\gcd(a,m)=1$ and $q\nmid m$
for every exceptional modulus
$q$ larger than $(\log x)^A$.
Both $C>0$ and the implied constant depend on $A$. The
dependence is ineffective. 

Zero-density results on $L$-functions, together with the standard zero-free
regions and bounds on $|\log L(s,\chi)|$, yield rather general results on the
sum of $\Lambda$, $\mu$, $\lambda$ and the like over arithmetic progressions
in short intervals. If we take as inputs the zero-density results in
\cite{Hux}, (1.1), and \cite{Pin}, Theorem 3.2, and the zero-free
region in \cite{Iw}, Theorem 2, then the procedure in \cite{Ra}
yields
\begin{equation}\label{eq:huxram}
\left|\mathop{\sum_{x < n\leq x+h}}_{n\equiv a \mo m} \lambda(n)\right| \ll
h e^{- C \frac{(\log x)^{1/4}}{(\log \log x)^{3/4}}} + x^{\frac{7}{12} + \epsilon}
\end{equation}
for any $x\geq 1$ and any $a$, $m$ such that $\gcd(a,m)=1$ and $q\nmid m$
for every exceptional modulus
$q$ larger than $(\log x)^A$.
Both $C>0$ and the implied constant depend on $A$ and $\epsilon$. The
dependence is ineffective. 

Neither (\ref{eq:sw}) nor (\ref{eq:huxram}) 
give the best upper bounds that could conceivably be given nowadays;
the zero-free region in \cite{Iw} can be used to improve the
exponent $-C \sqrt{\log x}$ in (\ref{eq:sw}) to 
$- C \frac{(\log x)^{4/7}}{(\log \log x)^{3/7}}$, provided that
the modulus $m$ is $\leq (\log x)^{3/7} (\log \log x)^{3/7}$. 
In fact, it seems reasonable to hope that, given current technology,
one might be able to sharpen the result in \cite{Iw}
 to give a zero-free region
for Dirichlet $L$-functions
almost as broad as the Korobov-Vinogradov region for the zeta function,
though there is a gap in the literature at this point. If we could
prove the existence of such a zero-free region,
 we could improve the exponent in (\ref{eq:sw}) 
to $- C (\log x)^{3/5 - \epsilon}$ and the exponent in (\ref{eq:huxram})
to $- C (\log x)^{1/3 - \epsilon}$.
Both (\ref{eq:sw}) and (\ref{eq:huxram}), however,
are already
stronger than what we need; in particular, Ingham's $\frac{5}{8}$ or Hoheisel's
$1-\delta$, $\delta>0$, would do in place of Huxley's $\frac{7}{12}$.
We will certainly do without any hypothetical refinements of existing results.

The proof of Theorem 3.3 does not require the distinction between
exceptional and non-exceptional moduli to be made in the application of
(\ref{eq:huxram}), since the modulus $q$ will always be smaller than
$(\log x)^A$ for some constant $A$.
Short-interval estimates will not be needed in the proofs of 
Theorem \ref{thm:two} or Proposition \ref{prop:three}.
One may avoid the use of short-interval
estimates altogether by setting $M_2 = e^{\epsilon \sqrt{\log N}}$ 
in the proof of Theorem \ref{thm:one} and proceeding as in the
proof of Proposition \ref{prop:three}, at the cost of replacing the
factor of $\frac{\log \log N}{\log N}$ in the statement of the theorem
by a factor of $\frac{\log \log N}{\sqrt{\log N}}$; alternatively, one
may prove Theorem \ref{thm:one} as we prove Theorem \ref{thm:two}, 
namely, via the Bombieri-Vinogradov theorem, and thus keep the bound we
now have in the statement of Theorem \ref{thm:one}.

\subsection{Bilinear bounds}
We shall need bounds for bilinear sums involving
the Liouville function. For sections
\ref{sec:three} and \ref{sec:theman}, the following lemma will suffice. It is simply a linear
bound in disguise.
\begin{lem}\label{lem:diverto}
Let $S$ be a convex subset of $\lbrack -N, N\rbrack^2$. Let 
$L\subset \mathbb{Z}^2$ be a lattice coset of index 
$\lbrack \mathbb{Z}^2 : L\rbrack$ not divisible by any exceptional
moduli larger than $(\log N)^A$.
Let $f:\mathbb{Z}\to \mathbb{C}$ be a function with $\max_y |f(y)|\leq 1$.
 Then, for every $\epsilon>0$,
\begin{equation}\label{eq:gen}
\left| \sum_{(x,y)\in S\cap L} \lambda(x) f(y) \right| \ll
\Area(S) \cdot
 e^{- C \frac{(\log N)^{1/4}}{(\log \log N)^{3/4}}} + N^{\frac{19}{12}+\epsilon},\end{equation}
where $C$ and the implied constant in (\ref{eq:gen})
depend only on $A$ and $\epsilon$.
\end{lem}
\begin{proof}
For every $y\in \mathbb{Z}\cap \lbrack - N, N\rbrack$, the set 
$\{x : (x,y)\in L\}$ is either the empty set or an arithmetic progression
$m_y \mathbb{Z} + a_y$, where $m_y | \lbrack \mathbb{Z}^2 : L \rbrack$.
Let $y_0$ and $y_1$ 
be the least and the greatest $y\in \mathbb{Z}\cap \lbrack -N,N\rbrack$ 
such that $\{x : (x,y)\in S\}$ is non-empty. Let $y\in \mathbb{Z} \cap
\lbrack y_0,y_1\rbrack$. 
Since $S$ is convex and a subset of $\lbrack - N, N\rbrack^2$, 
the set
$\{x : (x,y) \in S\}$ is an interval $\lbrack N_{y,0},N_{y,1}\rbrack$ contained
in $\lbrack - N, N\rbrack$. Hence
\[\begin{aligned}
\left| \sum_{(x,y)\in S\cap L} \lambda(x) f(y) \right| &=
\left| \mathop{\sum_{y_0\leq y\leq y_1}}_{\{x : (x,y)\in L\} \ne \emptyset}
\;\;\mathop{\sum_{N_{y,0}\leq x\leq N_{y,1}}}_{x\equiv a_y \mo m_y}
 \lambda(x) f(y) \right| \\ &\leq
\mathop{\sum_{y_0\leq y\leq y_1}}_{\{x : (x,y)\in L\} \ne \emptyset}
\left| \mathop{\sum_{N_{y,0}\leq x\leq N_{y,1}}}_{x\equiv a_y \mo m_y}
 \lambda(x) \right| .\end{aligned}\]
By (\ref{eq:sw}) and (\ref{eq:huxram}),
\[\begin{aligned}
 \mathop{\sum_{y_0\leq y\leq y_1}}_{\{x : (x,y)\in L\} \ne \emptyset}
\left| \mathop{\sum_{N_{y,0}\leq x\leq N_{y,1}}}_{x\equiv a_y \mo m_y}
 \lambda(x) \right| 
 &\ll \sum_{y_0\leq y\leq y_1} (N_{y_1}-N_{y_0}) e^{- C \frac{(\log N)^{1/4}}{
(\log \log N)^{3/4}}}
 + N^{\frac{19}{12}+\epsilon}
\end{aligned}\]
for any $\epsilon>0$.
Clearly \[\Area(S) = \sum_{y=y_0}^{y_1} (N_{y,1} - N_{y,0}) + O(N) .\]
Therefore
\[
\left| \sum_{(x,y)\in S\cap L} \lambda(x) f(y) \right| \ll
\Area(S) \cdot e^{- C \frac{(\log N)^{1/4}}{(\log \log N)^{3/4}}} +
N^{\frac{19}{12} +\epsilon} .\]
\end{proof}
As a special case of, say, Theorem 1 in \cite{Le}, we have the following
analogue of Bombieri-Vinogradov:
\begin{equation}\label{eq:bv}
\sum_{m\leq \frac{N^{1/2}}{(\log N)^{2 A+4}}} \mathop{\max_a}_{(a,m)=1} \max_{x\leq N}
\left|\mathop{\sum_{n\leq x}}_{n\equiv a \mo m} \lambda(n) -
\frac{1}{\phi(m)} \mathop{\sum_{n\leq x}}_{\gcd(n,m)=1} \lambda(n)\right| 
\ll \frac{N}{(\log N)^A},\end{equation}
where the implied constant depends only on $A$.

A simpler statement is true.
\begin{lem}\label{lem:bv2}
For any $A>0$,
\[\sum_{m\leq \frac{N^{1/2}}{(\log N)^{2 A+6}}} \max_a 
\max_{x\leq N} \left|\mathop{\sum_{n\leq x}}_{n\equiv a \mo m} \lambda(n) \right| 
\ll \frac{N}{(\log N)^A},\] where the implied constant depends only on $A$.
\end{lem}
\begin{proof}
Write $\rad(m) = \prod_{p|m} p$. Then
\[\sum_{d|\gcd(\rad(m),n)} \lambda(n/d) = \begin{cases}
\lambda(n) & \text{if $\gcd(m,n)=1$}\\
0 &\text{otherwise.}\end{cases}\]
Therefore
\[\begin{aligned}
\sum_{m\leq N^{1/2}} \frac{1}{\phi(m)}
 \max_{x\leq N} \left|\mathop{\sum_{n\leq x}}_{\gcd(n,m)=1} \lambda(n) 
\right| &=
\sum_{m\leq N^{1/2}} \frac{1}{\phi(m)}
 \max_{x\leq N} \left|\sum_{d|\rad(m)} 
\mathop{\sum_{n\leq x}}_{d|n} \lambda(n/d)\right|\\
&\leq \sum_{m\leq N^{1/2}} \frac{1}{\phi(m)}
 \sum_{d|\rad(m)} \max_{x\leq N/d} \left|\sum_{n\leq x} \lambda(n)\right|
\\ &\ll \sum_{m\leq N^{1/2}} \frac{1}{\phi(m)}
 \sum_{d|\rad(m)} N/d \cdot e^{-C \sqrt{\log N/d}}\\
&\leq N e^{- C \sqrt{\log N^{1/2}}} \sum_{m\leq N^{1/2}} \frac{1}{\phi(m)}
\sum_{d|\rad(m)} \frac{1}{d} 
\\ &\ll \frac{N}{(\log N)^A}.\end{aligned}\]
By (\ref{eq:bv}) this implies
\[\sum_{m\leq \frac{N^{1/2}}{(\log N)^{2 A+6}}} \mathop{\max_a}_{\gcd(a,m)=1} 
\max_{x\leq N} \left|\mathop{\sum_{n\leq x}}_{n\equiv a \mo m} \lambda(n) \right| 
\ll \frac{N}{(\log N)^A} .\]
Now
\[\begin{aligned}
\sum_{m\leq \frac{N^{1/2}}{(\log N)^{2 A+6}}} &\max_a 
\max_{x\leq N} \left|\mathop{\sum_{n\leq x}}_{n\equiv a \mo m} \lambda(n) \right| \\ &=
\sum_{m\leq \frac{N^{1/2}}{(\log N)^{2 A+6}}} \max_{r|m}
\max_{(a,m)=1} 
\max_{x\leq N} \left|\mathop{\sum_{n\leq x}}_{n\equiv a r\mo m} \lambda(n) \right| \\
%&=
%\sum_{m\leq \frac{N^{1/2}}{(\log N)^{2 A+6}}} \max_{r|m}
%\max_{(a,m)=1} 
%\max_{x\leq \frac{N}{r}} \left|\mathop{\sum_{n\leq x}}_{n\equiv a\mo m/r} \lambda(n) \right|\\
&< \sum_{r\leq N^{1/2}}
\sum_{s\leq \frac{(N/r)^{1/2}}{(\log (N/r))^{2 A+6}}}
\max_{(a,s)=1}
\max_{x\leq \frac{N}{r}} \left|\mathop{\sum_{n\leq x}}_{n\equiv a\mo s} \lambda(n) \right|
\\ &\ll \sum_{r\leq N^{1/2}} \frac{N/r}{(\log N/r)^{A+1}} \ll
\frac{N}{(\log N)^A} .\end{aligned}\]
\end{proof}

The following lemma is to Lemma \ref{lem:diverto} what Bombieri-Vinogradov is to (\ref{eq:sw}).

\begin{lem}\label{bomb4}
Let $A$, $K$ and $N$ be positive and satisfy
$K\leq \sqrt{N}/(\log N)^{2 A + 6}$. For $j=1,2,\dotsc,K$,
let $S_j$ be a convex subset of $\lbrack -N,N\rbrack^2$ and let
$L_j\subset \mathbb{Z}^2$ be a lattice coset of index $j$. Let $f:\mathbb{Z}\to \mathbb{C}$ be
a function with $\max_y |f(y)| \leq 1$. Then
\[\sum_{j=1}^{K} \left|
\sum_{(x,y)\in S_j\cap L_j} \lambda(x) f(y) \right|
\ll \frac{N^2}{(\log N)^A},\]
where the implicit constant depends only on $A$.
\end{lem}
\begin{proof}
We start with
\[\begin{aligned}
\sum_{j=1}^{K} \left|
\sum_{(x,y)\in S_j \cap L_j} \lambda(x) f(y) \right| &\leq
\sum_{j=1}^{K} 
 \sum_y \left| \mathop{\sum_x}_{(x,y)\in S_j\cap L_j} \lambda(x) \right|\\
&= 
\sum_{j=1}^{K}
 \sum_{k=0}^{\lceil N/j \rceil} \sum_{y = k j}^{(k+1) j -1} \left|
\mathop{\sum_x}_{(x,y)\in S_j\cap L_j} \lambda(x) \right| .\end{aligned}\]

For any $y\in \mathbb{Z}$, the set 
$\{x : (x,y)\in L_j\}$
is either the empty set or
an arithmetic progression of modulus $m_j|j$ independent of $y$.
Thus the set
\[A_j=\{(x,y)\in L_j: k j\leq y\leq (k+1) j-1\}\] is the union of $m_j$ sets
of the form
\[B_{y_0,a} = \{(x,y)\in \mathbb{Z}^2 : x\equiv a \mo m_j,\, y = y_0\}\]
with $k j \leq y_0 \leq (k+1) j - 1$. Since an arithmetic progression
of modulus $d$ is the union of $j/d$ arithmetic progressions of modulus $j$, the 
set $A_j$ is the union of $j$ sets of the form
\[C_{x_0,a} =  \{(x,y)\in \mathbb{Z}^2 : x\equiv a \mo j,\, y = y_0\}.\]
Therefore
\[\begin{aligned}
\sum_{j=1}^K
 \sum_{k=0}^{\lceil N/j \rceil} \sum_{y = k j}^{(k+1) j -1} \left|
\mathop{\sum_x}_{(x,y)\in S_j\cap L_j} \lambda(x) \right| &\leq
\sum_{j=1}^K
 \sum_{k=0}^{\lceil N/j \rceil} \sum_{l = 1}^j \left|
\mathop{\sum_x}_{(x,y_0(k,l))\in S_j\cap C_{y_0(k,l),a(k,l)}} 
 \lambda(x) \right| \\ &\leq \sum_{j=1}^K
(N+j)
\max_{y_0} \max_a \left|
\mathop{\sum_x}_{(x,y_0)\in S\cap C_{y_0,a}} \lambda(x)\right|\\ &\leq
\sum_{j=1}^K
(N+j)
\max_{-N\leq b\leq c\leq N} \max_a \left|
\mathop{\sum_{b\leq x\leq c}}_{x\equiv a \mo j} \lambda(x)\right|\\
&\leq \sum_{j=1}^K
 4 (N+j)
\max_{0< c\leq N} \max_a \left|
\mathop{\sum_{0<x\leq c}}_{x\equiv a \mo j} \lambda(x)\right|.\end{aligned}\]

We apply Lemma \ref{lem:bv2} and are done.
\end{proof}

\begin{cor}\label{bomb5}
Let $A$, $K$, $N$, $d_0$ and $d_1$ be positive integers such that
$K d_1$ is no larger than $N^{1/2}/(\log N)^{2 A + 6}$. For $k=1,2,\dotsc ,K$,
let $S_k$ be a convex subset of $\lbrack -N,N\rbrack^2$ and let
$L_k\subset \mathbb{Z}^2$ be a lattice coset of index $\frac{r_k}{d_0} k$ for
some $r_k$ dividing $d_0 d_1$. Then
\[\sum_{k\leq K} \left|
\sum_{(x,y)\in S_k\cap L_k} \lambda(x) \lambda(y) \right|
\ll \tau(d_0 d_1)\cdot \frac{N^2}{(\log N)^A},\]
where the implicit constant depends only on $A$. 
\end{cor}
\begin{proof}
For every $j\leq K d_1$, there are at most $\tau(d_0 d_1)$ lattice cosets
$L_k$ of index $j$. There are no lattice cosets $R_k$ of index greater than
$K d_1$. The statement then follows
from Lemma \ref{bomb4}.
\end{proof}

\section{The average of $\lambda$ on the product of three linear 
factors}\label{sec:three}
\begin{lem}\label{split1}
Let $\mathscr{P}\subset \{M_1,M_1+1,\dotsc,M_2-1\}$ 
be a set of primes between the positive integers
$M_1$ and $M_2$. Then there are
  $\sigma_d\in \mathbb{R}$ with $|\sigma_d|\leq 1$
and support on
\[\{M_1\leq d < M_2 : p|d \Rightarrow p\in \mathscr{P}\}\]
such that
\[\begin{aligned}\sum_{(x,y)\in S\cap L} g(x) f(x,y) &=
  \mathop{\sum_a \sum_b \sum_c}_{(a b,c)\in S\cap L} \sigma_a\,
g(a) g(b) f(a b, c) \\ &+ O\left(
\mathop{\prod_{M_1\leq p< M_2}}_{p\notin \mathscr{P}}
\left( 1 - \frac{1}{p}\right)^{-1} \cdot
\frac{\log M_1}{\log M_2} \frac{\Area(S)}{
\lbrack \mathbb{Z}^2:L\rbrack}
  + N M_2\right) \end{aligned}\]
for any positive integer $N>M_2$,
any convex set $S\subset \lbrack -N, N\rbrack^2$,
any lattice coset $L\subset \mathbb{Z}^2$ 
with index
$\lbrack \mathbb{Z}^2 : L\rbrack < M_1$,
any function $f:\mathbb{Z}^2\to \mathbb{C}$ 
and any completely multiplicative function
$g:\mathbb{Z}^2 \to \mathbb{C}$ with 
\[\max_{x,y} |f(x,y)|\leq 1,\;\; \max_y |g(y)|\leq 1.\] The implied constant is absolute.
\end{lem}
\begin{proof}
Let $y_1 = \min(\{y\in \mathbb{Z}:\exists x \,\st\, (x,y)\in S\cap L\})$. 
There is an
$l|\lbrack \mathbb{Z}^2:L\rbrack$ such that, for any $y\in \mathbb{Z}$, 
\[(\exists x \,\st\, (x,y)\in L) \Leftrightarrow
(l|y-y_1) .\]
Let
\[\begin{aligned}
N_{j,0} &= \min(\{x: (x, y_1 + j l) \in S\cap L\})\\
N_{j,1} &= \max(\{x: (x, y_1 + j l) \in S\cap L\})+1 .
\end{aligned}\]
Now take $\sigma_d$ as in Lemma \ref{lem:sieve}. 
Then
\[
\sum_{x : (x, y_1 + j l) \in S\cap L}
\left|1 - \sum_{d|x} \sigma_d\right| \ll
\mathop{\prod_{M_1\leq p< M_2}}_{p\notin \mathscr{P}} 
\left(1 - \frac{1}{p}\right)^{-1} \cdot
\frac{\log M_1}{\log M_2} \frac{N_{j,1}-N_{j,0}}{
\lbrack \mathbb{Z}^2 : L\rbrack/l} + M_2\]
Summing this over all $j$ we obtain
\[\begin{aligned}
\sum_{(x,y)\in S\cap L}
\left|1 - \sum_{d|x} \sigma_d \right|&\ll
\mathop{\prod_{M_1\leq p< M_2}}_{p\notin \mathscr{P}} 
\left(1 - \frac{1}{p}\right)^{-1} \cdot
\frac{\log M_1}{\log M_2} \frac{(\Area(S))/ l}{\lbrack \mathbb{Z}^2 : L
\rbrack/l} + M_2 N \\ &\ll 
\mathop{\prod_{M_1\leq p< M_2}}_{p\notin \mathscr{P}} 
\left(1 - \frac{1}{p}\right)^{-1} \cdot
\frac{\log M_1}{\log M_2} \frac{\Area(S)}{\lbrack \mathbb{Z}^2 : L \rbrack} 
+ M_2 N .
\end{aligned}\]
Since
\[
\left|\sum_{(x,y)\in S\cap L} g(x) f(x,y) -
\sum_{(x,y)\in S\cap L} \sum_{d|x} \sigma_d g(x) f(x,y) \right|\]
is at most
\[\sum_{(x,y)\in S\cap L} \left| g(x) f(x,y) -
\sum_{d|x} \sigma_d g(x) f(x,y) \right|
\leq \sum_{(x,y)\in S\cap L} \left|1 - \sum_{d|x} \sigma_d\right|
\]
and
\[\mathop{\sum_a \sum_b \sum_c}_{(a b,c)\in S\cap L} \sigma_a\,
g(a) g(b) f(a b, c) = 
\sum_{(x,y)\in S\cap L} \sum_{d|x} \sigma_d g(x) f(x,y),\]
we are done.  
\end{proof}
\begin{lem}\label{mat3}
Let $c_1$, $c_2$ be integers. 
Let $L\subset \mathbb{Z}^2$ be a lattice. Then the set
$\{(a,b)\in \mathbb{Z}^2 : (a,b c_1),(a,b c_2)\in L\}$ 
is 
either the empty set or a lattice coset $L'\subset \mathbb{Z}^2$
of index dividing $\lbrack \mathbb{Z}^2:L\rbrack^2$.
\end{lem}
\begin{proof}
The set of all elements of $L$ of the form $(a,b c_1)$ is the intersection
of a lattice coset of index $\lbrack \mathbb{Z}^2:L\rbrack$ and a lattice of index $c_1$.
By (\ref{eq:intersl}) it
 is either the empty set or
a lattice coset of index dividing $c_1 \lbrack \mathbb{Z}^2:L\rbrack$. Therefore
the set of all $(a,b)$ such that $(a,b c_1)$ is in $L$ is either
the empty set or a lattice coset
$L_1$
of index dividing $\frac{1}{c_1} c_1 \lbrack \mathbb{Z}^2:L\rbrack = \lbrack \mathbb{Z}^2:L\rbrack$.
Similarly, the set of all $(a,b)$ such that $(a,b c_2)\in L$ is either the
empty set or a lattice coset
$L_2$ of index dividing $\lbrack \mathbb{Z}^2:L\rbrack$. Therefore $L'=L_1\cap L_2$
is either the empty set or 
a lattice coset of index dividing $\lbrack \mathbb{Z}^2:L\rbrack^2$. 
\end{proof}
\begin{defn}
For $A = \left(\begin{matrix} a_{1 1} & a_{1 2}\\ a_{2 1} & a_{2 2}\\a_{3 1} & a_{3 2}\end{matrix}\right)$ we denote
\[A_{1 2} = \left(\begin{matrix} a_{1 1} & a_{1 2}\\ a_{2 1} & a_{2 2}\end{matrix}\right)\;\;\;
A_{1 3} = \left(\begin{matrix} a_{1 1} & a_{1 2}\\ a_{3 1} & a_{3 2}\end{matrix}\right)\;\;\;
A_{2 3} = \left(\begin{matrix} a_{2 1} & a_{2 2}\\ a_{3 1} & a_{3 2}\end{matrix}\right) .\]
\end{defn}
\begin{thm}\label{thm:one}
Let $S$ be a convex subset of $\lbrack -N,N\rbrack^2$, $N>1$. Let $L\subset \mathbb{Z}^2$ be a lattice coset. Let $a_{1 1}$, $a_{1 2}$, $a_{2 1}$, $a_{2 2}$, $a_{3 1}$, $a_{3 2}$ be rational integers. Then
\[\sum_{(x,y)\in S\cap L} 
\lambda((a_{1 1} x + a_{1 2} y)(a_{2 1} x + a_{2 2} y)(a_{3 1} x + a_{3 2} y))
\ll \frac{\log \log N}{\log N} \frac{\Area(S)}{
\lbrack \mathbb{Z}^2 : L\rbrack}
+ \frac{N^2}{(\log N)^\alpha}\]
for any $\alpha>0$. The implied constant depends only on $(a_{i j})$
and $\alpha$.
\end{thm}
\begin{proof}
We can assume that $A_{12}$ is non-singular, as otherwise the
statement follows immediately from Lemma \ref{lem:diverto}. 
Changing variables we obtain
\[\begin{aligned}
\sum_{(x,y)\in S\cap L}
&\lambda(a_{1 1} x + a_{1 2} y) \lambda(a_{2 1} x + a_{2 2} y)
\lambda(a_{3 1} x + a_{3 2} y) \\&=
\sum_{(x,y)\in A_{1 2} S \cap A_{1 2} L}
\lambda(x) \lambda(y) \lambda\left((a_{3 1}\, a_{3 2}) A_{1 2}^{-1} 
\left(\begin{matrix}x\\y\end{matrix}\right)\right)\\&=
\sum_{(x,y)\in A_{1 2} S \cap A_{1 2} L}
\lambda(x) \lambda(y) \lambda(q_1 x + q_2 y) ,\end{aligned}\]
where $q_1 = -\frac{\det(A_{2 3})}{\det(A_{1 2})}$ and
$q_2 = \frac{\det(A_{1 3})}{\det(A_{1 2})}$. Note that
$q_1 x + q_2 y$ is an integer for all $(x,y)$ in $A_{1 2} L$. We can assume that neither $q_1$ nor $q_2$ is zero. Write $S' = A_{1 2} S$, $L' = A_{1 2} L$.
Clearly $S'\subset \lbrack -N',N'\rbrack^2$ for
  $N' = \max(|a_{1 1}| + |a_{1 2}|,|a_{2 1}|+|a_{2 2}|) N$.

Now let
\[M_1 = (\log N')^{2 \alpha + 2}, M_2 = (N')^{\beta},\]
where $\alpha>0$ and $\beta\in (0,1/2)$ will be set later. Clearly $M_2>M_1$
for $N>N_0$, where $N_0$ depends only on $(a_{ij})$, $\alpha$ and $\beta$.

 By Lemma \ref{split1} with $\mathscr{P} = \{M_1\leq p<M_2 :
\text{$p$ prime}\}$,
\[\begin{aligned}\sum_{(x,y)\in S'\cap L'} \lambda(x) \lambda(y) 
\lambda(q_1 x + q_2 y) &=
  \mathop{\sum_a \sum_b \sum_c}_{(a b,c)\in S'\cap L'} \sigma_a\,
\lambda(a) \lambda(b) \lambda(c) 
\lambda(q_1 a b + q_2 c)\\ &+ O\left(\frac{\log M_1}{\log M_2} \frac{\Area(S')}{\lbrack \mathbb{Z}^2:L'\rbrack} + N' M_2\right) .\end{aligned}\]

We need to split the domain:
\[  \mathop{\sum_a \sum_b \sum_c}_{(a b,c)\in S'\cap L'} \sigma_a\,
\lambda(a) \lambda(b) \lambda(c) 
\lambda(q_1 a b + q_2 c)
= \sum_{s=1}^{\lceil M_2/M_1 \rceil} T_s
,\]
where
\[T_s = \mathop{\sum_{a = s M_1}^{(s+1) M_1 -1}
        \sum_{|b|\leq N'/s M_1} \sum_c}_{(a b,c)\in S' \cap L'}
 \sigma_a\,
\lambda(a) \lambda(b) \lambda(c) 
\lambda(q_1 a b + q_2 c) .\]
By Cauchy's inequality, 
\[
T_s^2 \leq \frac{(N')^2}{s M_1}
 \sum_c \sum_{|b|\leq N'/s M_1} 
\left(\mathop{\sum_{s M_1\leq a<(s+1) M_1}}_{(a b, c)\in S'\cap L'}
\sigma_a \lambda(a) \lambda(q_1 a b + q_2 c)\right)^2.\]
Expanding the square and changing the order of summation, we get
\[\frac{(N')^2}{s M_1}
\sum_{s M_1 \leq a_1, a_2 < (s+1) M_1}
%sum_{a_1 = s M_1}^{(s+1) M_1 -1}
%sum_{a_2 = s M_1}^{(s+1) M_1 -1}
\sigma_{a_1} \sigma_{a_2}
\lambda(a_1 a_2) \mathop{\sum_c 
\sum_{|b|\leq N'/s M_1}}_{(a_i b,c)\in S'\cap L'} \lambda(q_1 a_1 b + q_2 c)
 \lambda(q_1 a_2 b + q_2 c) .\]
There are at most $M_1\cdot 2 N' \frac{N'}{s M_1}$ terms with
$a_1=a_2$. They contribute at most $\frac{2 (N')^4}{s^2 M_1}$
to $T_s^2$, and thus no more than $((N')^2/\sqrt{M_1}) \log M_2$ to the
sum $\sum_{s=1}^{\lceil M_2/M_1\rceil} T_s$.
It remains to bound
\[\mathop{\sum_{s M_1 \leq a_1, a_2 < (s+1) M_1}}_{a_1\ne a_2}
\sigma_{a_1} \sigma_{a_2}
\lambda(a_1 a_2) \mathop{\sum_c 
\sum_{|b|\leq N'/s M_1}}_{(a_i b,c)\in S'\cap L'} \lambda(q_1 a_1 b + q_2 c)
 \lambda(q_1 a_2 b + q_2 c) .\]

Since $|\sigma_a|\leq 1$ for all $a$,
the absolute value of this is at most
\[ 
\mathop{\sum_{a_1= s M_1}^{(s+1) M_1 -1}
\sum_{a_2= s M_1}^{(s+1) M_1 -1}}_{a_1\ne a_2} 
 \left|
\mathop{\sum_c \sum_b}_{(a_i b,c)\in S'\cap L'} \lambda(q_1 a_1 b + q_2 c)
 \lambda(q_1 a_2 b + q_2 c) \right|.\]

By Lemma \ref{mat3} we may write
$\{(b,c)\in \mathbb{Z}^2: (a_1 b, c),(a_2 b, c)\in S'\cap L'\}$ as
$S_{a_1,a_2}''\cap L_{a_1,a_2}''$ with $S_{a_1,a_2}''$ a convex subset of
$\lbrack -N'/\max(a_1,a_2),N'/\max(a_1,a_2)\rbrack \times
\lbrack -N',N'\rbrack$ and $L_{a_1,a_2}''\subset \mathbb{Z}^2$ a lattice coset
of index dividing $\lbrack \mathbb{Z}^2:L'\rbrack^2$. 
Hence we have the sum
\begin{equation}\label{eq:coal1}
\mathop{\sum_{a_1= s M_1}^{(s+1) M_1 -1}
\sum_{a_2= s M_1}^{(s+1) M_1 -1}}_{a_1\ne a_2} 
\left|\sum_{(b,c)\in S_{a_1,a_2}''\cap L_{a_1,a_2}''} 
\lambda(q_1 a_1 b + q_2 c)
 \lambda(q_1 a_2 b + q_2 c) \right|.\end{equation}

Set $S_{a_1,a_2} = 
  \left(\begin{matrix}q_1 a_1 &q_2\\q_1 a_2 &q_2\end{matrix}\right) S_{a_1,a_2}''$,
$L_{a_1,a_2} = 
  \left(\begin{matrix}q_1 a_1 &q_2\\q_1 a_2 &q_2\end{matrix}\right)
L_{a_1,a_2}''$, $N'' = (|q_1|+|q_2|) N'$.  Clearly $S_{a_1,a_2}$ is a convex subset of $\lbrack -N'',N''\rbrack^2$ with 
\[\Area(S_{a_1,a_2}) = |q_1 q_2 (a_1-a_2)|
\Area(S'') \leq |q_1 q_2| M_1 \frac{4 (N')^2}{s M_1} \ll \frac{N^2}{s},\]
 whereas $L_{a_1,a_2}\subset \mathbb{Z}^2$ is a lattice coset of index 
$|q_1 q_2 (a_1 - a_2)| \lbrack \mathbb{Z}^2 : L_{a_1,a_2}''\rbrack$. (That $L_{a_1,a_2}$ is inside $\mathbb{Z}^2$ follows from
our earlier remark that $q_1 x + q_2 y$ is an integer for all
$(x,y)$ in $A_{1 2} L$.) Now we have
\[
\mathop{\sum_{a_1= s M_1}^{(s+1) M_1 -1}
\sum_{a_2= s M_1}^{(s+1) M_1 -1}}_{a_1\ne a_2} 
 \left| 
\sum_{(v,w)\in S_{a_1,a_2}\cap L_{a_1,a_2}} 
\lambda(v)
 \lambda(w) \right|.\]

This is at most
\begin{equation}\label{eq:coal2}
M_1^2 \max_{s M_1\leq a < (s+1) M_1} 
    \mathop{\max_{-M_1\leq d\leq M_1}}_{d\ne 0}
 \left|
\sum_{(v,w)\in S_{a,a+d}\cap L_{a,a+d}} 
\lambda(v)
 \lambda(w) \right|.\end{equation}
We can assume that $\lbrack \mathbb{Z}^2 : L\rbrack < (\log N)^{\alpha}$,
as otherwise the bound we are attempting to prove is trivial. Hence
$\lbrack \mathbb{Z}^2 : L ''\rbrack \ll (\log N)^{2 \alpha}$.
By Lemma \ref{lem:diverto},
\begin{equation}\label{eq:coal3} \left|
\sum_{(v,w)\in S_{a,a+d}\cap L_{a,a+d}} 
\lambda(v)
 \lambda(w) \right|\ll
\frac{N^2}{s} \cdot e^{-\frac{1}{2} \frac{C (\log N'')^{1/4}}{
(\log \log N'')^{3/4}}} + (N'')^{19/12+\epsilon} 
.\end{equation}
It is time to collect all terms. The total is at most a constant times
\[\begin{aligned}
\frac{\log M_1}{\log M_2} \frac{\Area(S')}{\lbrack \mathbb{Z}^2 :
L'\rbrack} &+ N' M_2 +
\frac{(N')^2}{\sqrt{M_1}} \log M_2 \\ 
&+ N' N'' \sqrt{M_1}
\log M_2 \cdot e^{-\frac{1}{2} \frac{C (\log N'')^{1/4}}{
(\log \log N'')^{3/4}}}
+ N' (N'')^{19/24 + \epsilon} \sqrt{M_2}
,
\end{aligned}\]
where the constant depends only on $(a_{i j})$ and $\alpha$. 
Set $\beta = \frac{1}{3}$. Simplifying
we obtain
\[O\left(\frac{\log \log N}{\log N} \frac{\Area(S)}{\lbrack \mathbb{Z}^2 : L\rbrack}
+ \frac{N^2}{(\log N)^\alpha}\right).\]
\end{proof}
\section{The average of $\lambda$ on the 
product of a linear and a quadratic factor}\label{sec:threetwo}
We will be working with quadratic extensions $K/\mathbb{Q}$. 
We define
\[\begin{aligned}
 \jmath(x + y \sqrt{d}) = (x,y) \;\; &\text{if $d\not\equiv 1 \mo 4$,}\\
 \jmath(x + y \sqrt{d}) = (x-y, 2y) \; &\text{if $d\equiv 1 \mo 4$,}
\end{aligned}\]
where $x, y\in \mathbb{Q}$. 

For every $z\in \jmath^{-1}(\lbrack -N,N \rbrack^2)$,
\begin{equation}
|N_{K/\mathbb{Q}} z| \ll N^2 ,
\end{equation}
where the implied constant depends only on $K$. In general there is no 
implication in the opposite sense, as the norm need not be positive definite.
For $K = \mathbb{Q}(\sqrt{d})$, $d<0$,
\begin{equation}\label{eq:boundn1}
\#\{z\in \mathscr{O}_K : N_{K/\mathbb{Q}}(z) \leq A\} \ll A .
\end{equation}
For $K = \mathbb{Q}(\sqrt{d})$, $d>1$,
\begin{equation}\label{eq:boundn2}
\#\{z\in \jmath^{-1}(\lbrack -N,N \rbrack^2) : |N_{K/\mathbb{Q}}(z)| \leq A\}
\ll A (1 + \log N) .
\end{equation}
In either case the implied constant depends only on $d$.

\begin{lem}\label{lem:divbyno}
Let $\mathfrak{a}$ be an ideal in $\mathbb{Q}(\sqrt{d})/\mathbb{Q}$ divisible by no
rational integer $n>1$. Then for any positive $N$,
$y_0\in \lbrack -N,N \rbrack$,
\[\#\{(x,y_0)\in \lbrack -N,N\rbrack^2 : \jmath^{-1}(x,y_0)\in \mathfrak{a}\}
\leq 2 \lceil N/N_{K/\mathbb{Q}}(\mathfrak{a})\rceil .\]
\end{lem}
\begin{proof}
For every rational integer $r\in \mathfrak{a}$, $N_{K/\mathbb{Q}} \mathfrak{a}
| r$. Hence
\[\{x : \jmath^{-1}(x,y_0) \in \mathfrak{a}\}\]
is an arithmetic progression of modulus $N_{K/\mathbb{Q}} \mathfrak{a}$.
\end{proof}
\begin{thm}\label{thm:two}
Let $S$ be a convex subset of $\lbrack -N,N\rbrack^2$, $N>1$. Let $L\subset \mathbb{Z}^2$ be a lattice
coset. Let $a_1$, $a_2$, $a_3$, $a_4$, $a_5$ be rational integers
such that $a_1 x^2 + a_2 x y + a_3 y^2$ is irreducible. Then
\[\sum_{(x,y)\in S\cap L} \lambda((a_1 x^2 + a_2 x y + a_3 y^2) (a_4 x + a_5 y)) \ll
\frac{\log \log N}{\log N} \frac{\Area(S)}{\lbrack \mathbb{Z}^2 : L\rbrack} +
\frac{N^2}{(\log N)^\alpha}\]
for any $\alpha>0$. The implied constant depends only on $(a_{i j})$ and $\alpha$.
\end{thm}
\begin{proof}
Write $d$ for $a_1^2 - 4 a_0 a_2$,
$K/\mathbb{Q}$ for $\mathbb{Q}(\sqrt{d})/\mathbb{Q}$,
$\mathfrak{N} z$ for $N_{K/\mathbb{Q}} z$ and
$\overline{r + s\sqrt{d}}$ for
$r - s\sqrt{d}$.
By Lemma \ref{lem:quad} there are 
$\alpha_1, \alpha_2 \in
\mathscr{O}_K$ linearly independent over 
$\mathbb{Q}$ and a non-zero rational number $k$ such that
\[a_1 x^2 + a_2 x y + a_3 y^2 = k \mathfrak{N}(x \alpha_1 + y \alpha_2)
= k (x \alpha_1 + y \alpha_2) \overline{(x \alpha_1 + y\alpha_2)} .\]
Hence
\[ \sum_{(x,y)\in S\cap L} \lambda((a_1 x^2 + a_2 x y + a_3 y^2) (a_4 x + a_5 y))\] equals
\[\lambda(k) \sum_{(x,y)\in S\cap L} \lambda((x \alpha_1 + y\alpha_2)
\overline{(x \alpha_1 + y\alpha_2)} (a_4 x + a_5 y)) .\]

We
write $\Re (r+s \sqrt{d})$ for $r$, $\Im (r + s \sqrt{d})$ for
$s$.
Let $C = \left(\begin{matrix}\Re \alpha_1 &\Re \alpha_2\\
\Im \alpha_1&\Im \alpha_2\end{matrix}
\right)^{-1}$. Then $a_4 x + a_5 y = q z + \overline{q z}$
for $z = x \alpha_1 + y \alpha_2$,
\[ q = \frac{1}{2} (a_4 c_{1 1} + a_5 c_{2 1} + \frac{1}{\sqrt{d}} (a_4 c_{1 2} + a_5 c_{2 2})).\] 
Define $\phi_Q:\mathbb{Z}^2\to \mathscr{O}_K$ to be the mapping
$(x,y) \mapsto (x \alpha_1 + y \alpha_2)$. Let 
$L' = (\jmath \circ \phi_Q)(L)$. 
Since $\jmath \circ \phi_Q$ is linear, it extends to a map $\varphi:
\mathbb{R} \to \mathbb{R}$. Let $S' = \phi(S)$. Then
\[  \sum_{(x,y)\in S\cap L} \lambda((x \alpha_1 + y\alpha_2)
\overline{(x \alpha_1 + y\alpha_2)} (a_4 x + a_5 y)) =
\sum_{\jmath(z)\in S'\cap L'} 
 \lambda(z \overline{z} (q z + \overline{q z})) .\]
Note that $q z + \overline{q z}$ is an integer for all $z\in L'$.

Let $N'$ be the smallest integer greater than one such that
$S'\subset \lbrack -N',N'\rbrack^2$. (Note that $N'\leq c_1 N$, where
$c_1$ is a constant depending only on $\mathbb{Q}$.) Suppose $K/\mathbb{Q}$
is real. Then, by (\ref{eq:boundn2}),
\[
\#\{z\in \jmath^{-1}(\lbrack-N',N'\rbrack^2) : |\mathfrak{N} z| 
       \leq \frac{(N')^2}{(\log N)^{\alpha+1}}\}
\leq \frac{(N')^2 (1 + \log N') }{(\log N)^{\alpha + 1}} 
\ll
\frac{N^2}{(\log N)^\alpha}
.\]
The set
\[
\{P\in \lbrack -N',N'\rbrack^2 : |\mathfrak{N} (\jmath^{-1}(P))| 
  > \frac{(N')^2}{(\log N)^{\alpha+1}}\}
\]
is the region within a square and outside two hyperbolas. As such it is
the disjoint union of at most four convex sets. Hence the set
\[S'' = S' \cap \{P\in \lbrack - N', N'\rbrack^2 : 
|\mathfrak{N}(\jmath^{-1}(P))|
 > (N')^2/(\log N)^{\alpha+1} \}\]
is the disjoint union of at most four convex sets:
\[S'' = S_1 \cup S_2 \cup S_3 \cup S_4 .\]
In the following, $S^*$ will be $S_1$, $S_2$, $S_3$ or $S_4$, and as such
a convex set contained in $S''$.

Suppose now that $K/\mathbb{Q}$ is imaginary. Then the set
\[\{P\in \lbrack - N',N'\rbrack^2 : \mathfrak{N}(\jmath^{-1}(P))
 > (N')^2/(\log N)^{\alpha+1} \}\]
is the region within a square and outside the circle given by
\begin{equation}\label{eq:satan}
\{P : \mathfrak{N}(\jmath^{-1}(P)) = (N')^2/(\log N)^{\alpha+1} \}.
\end{equation}
We can circumscribe about (\ref{eq:satan}) a rhombus containing no more
than
\[O((N')^2/(\log N)^{\alpha+1})\]
integer points, where the implied constant depends only on $Q$. We then
quarter the region inside the square $\lbrack -N', N' \rbrack^2$ and outside
the rhombus, obtaining four convex sets $S_1$, $S_2$, $S_3$, $S_4$ inside
$S$. We let $S^*$ be 
$S_1$, $S_2$, $S_3$ or $S_4$.

For $K$ either real or imaginary, we now have a convex set 
$S^* \subset \lbrack -N', N'\rbrack^2$ such that, for any 
$z\in \mathscr{O}_K$, 
\[\jmath(z) \in S^* \Rightarrow \mathfrak{N} z > (N')^2/(\log N)^{\alpha+1} .\]
Our task is to bound
\[\mathop{\sum_{z\in \mathscr{O}_K}}_{\jmath(z) \in S^* \cap L'}
 \lambda(z \bar{z} (q z + \overline{q z})) .
\]

Set \[M_1 = (\log N)^{20(\alpha +1)},\;\; 
M_2 = \frac{N^{1/2}}{4 d |\num(\mathfrak{N} q)| 
 \lbrack \mathscr{O}_K : L' \rbrack^2
(\log N)^{18 \alpha +24}},\]
where, for a rational number $r$, $\num(r)$ stands for the numerator $a$
 of $r = a/b$, $\gcd(a,b)=1$.
 By Lemma \ref{sieve2} with $\mathscr{P} = \{\text{$\mathfrak{p}$ prime}:
M_1\leq N \mathfrak{p} < M_2\}$,
\begin{equation}\label{eq:satun}\begin{aligned}
\mathop{\sum_{z\in \mathscr{O}_K}}_{\jmath(z)\in S^*\cap L'} 
  \lambda(z \overline{z} (q z + \overline{q z})) &=
\mathop{\sum_{z\in \mathscr{O}_K}}_{\jmath(z)\in S^*\cap L'} 
\sum_{\mathfrak{d} | z}
\sigma_{\mathfrak{d}} \lambda(z \overline{z} (q z + \overline{q z})) \\
&+ O\left(\frac{\log M_1}{\log M_2}
\frac{\Area(S^*)}{\lbrack \mathscr{O}_K:L'\rbrack}
 + N' M_2^2 \right).\end{aligned}\end{equation}
 Let $N'' = (9/4 + |d|) (N')^2 $. Then 
$\jmath(z)\in \lbrack - N', N'\rbrack^2$
implies $|\mathfrak{N} z|\leq N''$. Since $\sigma_\mathfrak{d}=0$
when $N \mathfrak{d}<M_1$, the first term on the right of
(\ref{eq:satun}) equals
\[ \mathop{\sum_{\mathfrak{b}}}_{\mathfrak{N} \mathfrak{b}\leq N''/M_1} 
\lambda(\mathfrak{b} \overline{\mathfrak{b}})
   \mathop{\sum_{\mathfrak{a}}}_{\text{$\mathfrak{a} \mathfrak{b}$ principal}}
\sigma_{\mathfrak{a}} \lambda(\mathfrak{a} \overline{\mathfrak{a}})
\mathop{\sum_{(z) = \mathfrak{a} \mathfrak{b}}}_{\jmath(z)\in S^* \cap L'} 
 \lambda(q z + \overline{q z}) .\]
We need to split the domain:
\begin{equation}\label{eq:lomsuf}
\mathop{\sum_{\mathfrak{b}}}_{\mathfrak{N} \mathfrak{b}\leq N''/M_1} 
\lambda(\mathfrak{b} \overline{\mathfrak{b}})
   \mathop{\sum_{\mathfrak{a}}}_{\text{$\mathfrak{a} \mathfrak{b}$ principal}}
\sigma_{\mathfrak{a}} \lambda(\mathfrak{a} \overline{\mathfrak{a}})
\mathop{\sum_{(z) = \mathfrak{a} \mathfrak{b}}}_{\jmath(z)\in S^* \cap L'} 
 \lambda(q z + \overline{q z}) =
\sum_{s=1}^{\lceil \log_2(N''/M_1)\rceil} T_s,\end{equation}
where
\[T_s = \mathop{\sum_{\mathfrak{b}}}_{2^{s-1}\leq \mathfrak{N} 
\mathfrak{b}\leq 2^s} 
\lambda(\mathfrak{b} \overline{\mathfrak{b}})
   \mathop{\sum_{\mathfrak{a}}}_{\text{$\mathfrak{a} \mathfrak{b}$ principal}}
\sigma_{\mathfrak{a}} \lambda(\mathfrak{a} \overline{\mathfrak{a}})
\mathop{\sum_{(z) = \mathfrak{a} \mathfrak{b}}}_{\jmath(z)\in S^* \cap L'} 
 \lambda(q z + \overline{q z}) .\]
Notice that $\lambda(\mathfrak{b} \overline{\mathfrak{b}})$,
$\sigma_{\mathfrak{a}}$,  $\lambda(\mathfrak{a} \overline{\mathfrak{a}})$ and
$\lambda(q z + \overline{q z})$ are all real.
By Cauchy's inequality,
\[\begin{aligned}
T_s^2&\ll 2^{s-1} \mathop{\sum_{\mathfrak{b}}}_{2^{s-1}\leq \mathfrak{N}
 \mathfrak{b}\leq 
2^s}
\left(\mathop{\sum_{\mathfrak{a}}}_{\text{$\mathfrak{a} \mathfrak{b}$ principal}}
\sigma_{\mathfrak{a}} \lambda(\mathfrak{a} \overline{\mathfrak{a}})
\mathop{\sum_{(z) = \mathfrak{a} \mathfrak{b}}}_{\jmath(z)\in S^* \cap L'} 
 \lambda(q z + \overline{q z})\right)^2 \\
&\leq
2^{s-1} \sum_{\mathfrak{b}}
\left(
\mathop{\mathop{\sum_{\mathfrak{a}}}_{\text{$\mathfrak{a} \mathfrak{b}$ principal}}}_{
n_{s 0}< \mathfrak{N} \mathfrak{a} \leq n_{s 1}}
\sigma_{\mathfrak{a}} \lambda(\mathfrak{a} \overline{\mathfrak{a}})
\mathop{\sum_{(z) = \mathfrak{a} \mathfrak{b}}}_{\jmath(z)\in S^* \cap L'} 
 \lambda(q z + \overline{q z})\right)^2 ,
\end{aligned}\]
where $n_{s 0} = \frac{(N')^2}{2^s (\log N)^{\alpha+1}}$ and
$n_{s 1} = \min(\frac{N''}{2^{s-1}}, M_2)$.
Expanding the square and changing the order of summation, we get
\[\begin{aligned}
2^{s-1} \mathop{\sum_{\mathfrak{a}_1}}_{n_{s 0} < \mathfrak{N} \mathfrak{a}_1\leq n_{s 1}}
&\mathop{\sum_{\mathfrak{a}_2}}_{n_{s 0} < \mathfrak{N} \mathfrak{a}_2
\leq n_{s 1}}
\sigma_{\mathfrak{a}_1} \sigma_{\mathfrak{a}_2}
\lambda(\mathfrak{a}_1 \overline{\mathfrak{a}_1})
\lambda(\mathfrak{a}_2 \overline{\mathfrak{a}_2})
\\ &\mathop{\sum_{\mathfrak{b}} }_{
\text{$\mathfrak{a}_1 \mathfrak{b}, \mathfrak{a}_2 \mathfrak{b}$ principal}}
\mathop{\sum_{(z_1) = \mathfrak{a}_1 \mathfrak{b}}}_{\jmath(z_1)\in S^* 
\cap L'} 
\mathop{\sum_{(z_2) = \mathfrak{a}_2 \mathfrak{b}}}_{\jmath(z_2)\in S^* 
\cap L'} 
 \lambda(q z_1 + \overline{q z_1}) \lambda(q z_2 + \overline{q z_2}) .
\end{aligned}\]

Write $\mathscr{S}(x + y\sqrt{d})$ for $\max(|x|,|y|)$. Let
$r = (z_2/z_1)\cdot \mathfrak{N} \mathfrak{a}$. We have
$r\in \overline{\mathfrak{a}_1}$ because 
\[(r) = ((z_2)/(z_1))\cdot \mathfrak{N} \mathfrak{a}_1 =
(\mathfrak{a}_2 / \mathfrak{a}_1)\cdot \mathfrak{N} \mathfrak{a}_1 =
\mathfrak{a}_2 \cdot \overline{\mathfrak{a}_1} .\]
Since $\mathfrak{N} z_1 > \frac{(N')^2}{(\log N)^{\alpha + 1}}$
and $\mathscr{S} (z_2 \overline{z_1}) \ll (N')^2$, where the implied
constant depends only on $\mathbb{Q}$,
\begin{equation}\label{eq:S}
\mathscr{S} (r) =
\mathscr{S}\left(\frac{z_2}{z_1} \mathfrak{N} \mathfrak{a}_1\right) =
\mathscr{S}\left(\frac{z_2 \overline{z_1}}{\mathfrak{N} z_1} 
\mathfrak{N} \mathfrak{a}_1\right) =
\mathscr{S}(z_2 \overline{z_1}) \frac{\mathfrak{N} \mathfrak{a}}{
\mathfrak{N} z_1} \ll n_{s 1} (\log N)^{\alpha + 1} .\end{equation}
Set \[R_s = \jmath^{-1}\left(\left[- k n_{s 1} (\log N)^{\alpha + 1},
k n_{s 1} (\log N)^{\alpha + 1}\right]^2\right) ,\]
where $k$ is the implied constant in (\ref{eq:S}) and as such depends only on $K$.
Changing variables we obtain
\[\begin{aligned}
2^{s-1} \mathop{\sum_{\mathfrak{a}}}_{n_{s 0} < \mathfrak{N}
 \mathfrak{a}\leq n_{s 1}}
&\mathop{\sum_{r\in \overline{\mathfrak{a}} \cap R_s}}_{
 n_{s 0} < \mathfrak{N} \left(\frac{(r)}{\mathfrak{a}}\right) \leq n_{s 1}}
\sigma_{\mathfrak{a}} \sigma_{(r)/\mathfrak{a}}
\lambda(\mathfrak{a} \overline{\mathfrak{a}})
\lambda\left(\frac{(r)}{\mathfrak{a}}
\overline{\frac{(r)}{\mathfrak{a}}}\right) \\
&\mathop{\mathop{\sum_z}_{\jmath(z)\in \jmath(\mathfrak{a})\cap S^*\cap L'}}_{
\jmath(r z/\mathfrak{N} \mathfrak{a}) \in S^*\cap L'} \lambda(q z + \overline{q z}) \lambda\left(\frac{q r z}{\mathfrak{N} \mathfrak{a}} + \overline{\frac{q r z}{\mathfrak{N} \mathfrak{a}}}\right) ,
\end{aligned}\]
that is, $2^{s-1}$ times
\begin{equation}\label{eq:somename}
\mathop{\sum_{\mathfrak{a}}}_{n_{s 0} < \mathfrak{N} \mathfrak{a} \leq n_{s 1}}
\mathop{\sum_{r\in \overline{\mathfrak{a}} \cap R_s}}_{
 n_{s 0} < \mathfrak{N} \left(\frac{(r)}{\mathfrak{a}}\right) \leq n_{s 1}}
\sigma_{\mathfrak{a}} \sigma_{(r)/\mathfrak{a}}
\lambda(r \overline{r})
\mathop{\mathop{\sum_{z}}_{\jmath(z)\in \jmath(\mathfrak{a})\cap S^*\cap L'}}_{
\jmath(r z/\mathfrak{N} \mathfrak{a})\in S^*\cap L'} \lambda(q z + \overline{q z}) \lambda\left(\frac{q r z}{\mathfrak{N} \mathfrak{a}} + 
\overline{\frac{q r z}{\mathfrak{N} \mathfrak{a}}}\right) .
\end{equation}
We now wish to eliminate the terms coming from $\mathfrak{a}$
with non-trivial rational integer divisors; we may do so once
we show that the total contribution of such terms is small.
For any non-zero rational integer $n$,
\[\begin{aligned}
\mathop{\mathop{\sum_{\mathfrak{a}}}_{n_{s 0} <\mathfrak{N} \mathfrak{a}
\leq n_{s 1}}}_{n|\mathfrak{a}}
\sum_{r\in \overline{\mathfrak{a}} \cap R_s}\;
\mathop{\sum_{z}}_{\jmath(z)\in \jmath(\mathfrak{a})\cap S^*} 1 &\ll
\mathop{\mathop{\sum_{\mathfrak{a}}}_{n_{s 0} <\mathfrak{N} \mathfrak{a}
\leq n_{s 1}}}_{n|\mathfrak{a}}
\frac{(2 k n_{s 1} (\log N)^{\alpha+1})^2}{\mathfrak{N} \mathfrak{a}} 
\left(\frac{(N')^2}{\mathfrak{N} \mathfrak{a}} + N'\right)
 \\
&\ll \frac{1}{n^2} \frac{N^4 (\log N)^{3\alpha +3}}{2^s}.
\end{aligned}\]
(Here $(N')^2/\mathfrak{N} \mathfrak{a} \gg N'$ because $\mathfrak{N}
\mathfrak{a} \leq n_{s 1}$ and $n_{s 1} \leq M_2$.)
Since the support of $\sigma_{\mathfrak{d}}$ is a subset of
\[\{\mathfrak{d} : M_1\leq \mathfrak{N} \mathfrak{d} < M_2,\;\;\;
\mathfrak{N} \mathfrak{p} < M_1 \Rightarrow \mathfrak{p}
\nmid \mathfrak{d}\},\]
we have that $n|\mathfrak{a}$ and $\sigma_{\mathfrak{a}} \ne 0$
imply $n\geq \sqrt{M_1}$. Therefore (\ref{eq:somename}) equals
\begin{equation}\label{eq:somename2}
\mathop{
 \mathop{\sum_{\mathfrak{a}}}_{
n_{s 0} < \mathfrak{N} \mathfrak{a}\leq n_{s 1}}}_{n>1
\Rightarrow n\nmid \mathfrak{a}}
\mathop{\sum_{r\in \overline{\mathfrak{a}} \cap R_s}}_{
 n_{s 0} < \mathfrak{N} \left(\frac{(r)}{\mathfrak{a}}\right) \leq
 n_{s 1}}
\sigma_{\mathfrak{a}} \sigma_{(r)/\mathfrak{a}}
\lambda(r \overline{r})
\mathop{\mathop{\sum_{z}}_{\jmath(z)\in \jmath(\mathfrak{a})\cap S^*\cap L'}}_{
\jmath(r z/\mathfrak{N} \mathfrak{a})\in S^*\cap L'} \lambda(q z + \overline{q z}) \lambda\left(\frac{q r z}{\mathfrak{N} \mathfrak{a}} + \overline{\frac{q r z}{\mathfrak{N} \mathfrak{a}}}\right) 
\end{equation}
plus $O(N^4 (\log N)^{3 \alpha +3}/(2^s \sqrt{M_1}))$. The absolute value of
(\ref{eq:somename2}) is at most
\begin{equation}\label{eq:somename3}
\mathop{\mathop{\sum_{\mathfrak{a}}}_{
n_{s 0} < \mathfrak{N} \mathfrak{a} \leq n_{s 1}}}_{n>1\Rightarrow n\nmid \mathfrak{a}}
\sum_{r\in \overline{\mathfrak{a}} \cap R_s}
\left|\mathop{\mathop{\sum_{z}}_{\jmath(z)\in \jmath(\mathfrak{a})\cap S^*\cap L'}}_{
\jmath(r z/\mathfrak{N} \mathfrak{a})\in S^*\cap L'}
 \lambda(q z + \overline{q z}) \lambda\left(\frac{q r z}{\mathfrak{N} \mathfrak{a}} + \overline{\frac{q r z}{\mathfrak{N} \mathfrak{a}}}\right) \right| .
\end{equation}
By Lemma \ref{lem:divbyno},
\[\begin{aligned}
\mathop{\mathop{\sum_{\mathfrak{a}}}_{
n_{s 0} < \mathfrak{N} \mathfrak{a}
\leq n_{s 1}}}_{n>1 \Rightarrow n\nmid \mathfrak{a}}
\sum_{r\in \overline{\mathfrak{a}} \cap R_s\cap \mathbb{Z}}
\mathop{\sum_{z}}_{\jmath(z)\in \jmath(\mathfrak{a})\cap S^*} 1 \;\;&\ll
\mathop{\sum_{\mathfrak{a}}}_{
n_{s 0} < \mathfrak{N} \mathfrak{a}
\leq n_{s 1}} \frac{k n_{s 1} (\log N)^{\alpha+1}}{
\mathfrak{N} \mathfrak{a}} \cdot
 \left(\frac{(N')^2}{\mathfrak{N} \mathfrak{a}} + N'\right) \\ 
%&\ll
%\frac{N^4 (\log N)^{\alpha+1}}{2^{s-1} n_{s 0}} +
%\frac{N^3 (\log N)^{\alpha+2}}{2^{s-1}} \\ 
&\ll
N^2 (\log N)^{2(\alpha+1)} 
.\end{aligned}\]

Thus we are left with
\begin{equation}\label{eq:somename4}
\mathop{\mathop{\sum_{\mathfrak{a}}}_{n_{s 0} < \mathfrak{N} \mathfrak{a} \leq n_{s 1}}}_{n>1\Rightarrow n\nmid \mathfrak{a}}
\mathop{\sum_{r\in \overline{\mathfrak{a}} \cap R_s}}_{\Im r \ne 0}
\left|\mathop{\mathop{\sum_{z}}_{\jmath(z)\in \jmath(\mathfrak{a})\cap S^*\cap L'}}_{
\jmath(r z/\mathfrak{N} \mathfrak{a})\in S^*\cap L'}
 \lambda(q z + \overline{q z}) \lambda\left(\frac{q r z}{
\mathfrak{N} \mathfrak{a}} + 
\overline{\frac{q r z}{\mathfrak{N} \mathfrak{a}}}\right) \right| .
\end{equation}

Notice that $r\in \overline{\mathfrak{a}}$ and
$z\in \mathfrak{a}$ imply $(r z/\mathfrak{N} \mathfrak{a})\in \mathscr{O}_K$. Hence
$(r/\mathfrak{N} \mathfrak{a})^{-1} \mathscr{O}_K \supset \mathfrak{a}$. Therefore
$(r/\mathfrak{N} \mathfrak{a})^{-1} \jmath^{-1}(L') \cap \mathfrak{a}$ is either the empty set or a sublattice of 
$\mathfrak{a}$ of index dividing $\lbrack \mathscr{O}_K : L' \rbrack$. This means that
\[L_{\mathfrak{a},r} = \{z \in \mathfrak{a}\cap \jmath^{-1}(L'): (r z/\mathfrak{N} \mathfrak{a})\in \jmath^{-1}(L')\}\]
is either the empty set or a sublattice of $\mathfrak{a}$ of index 
$\lbrack \mathfrak{a} : L_{\mathfrak{a},r} \rbrack$ 
dividing $\lbrack \mathscr{O}_K : L' \rbrack^2$,
whereas \[S_{\mathfrak{a},r} = \{\jmath(z)
 \in S^* : \jmath(r z/\mathfrak{N} \mathfrak{a}) \in S^*\}\] 
is a convex subset of $\lbrack -N',N' \rbrack^2$.
The map 
\[\kappa:(x,y) \mapsto \left(
q \cdot \jmath^{-1}(x,y) + \overline{q \cdot \jmath^{-1}(x,y)},
\frac{q r \cdot \jmath^{-1}(x,y)}{\mathfrak{N} \mathfrak{a}} + \overline{\frac{q r 
\cdot \jmath^{-1}(x,y)}{\mathfrak{N} \mathfrak{a}}}\right)\]
is given by the matrix
\[\begin{aligned}
\left(\begin{matrix} 2 & 0\\ 2\frac{\Re r}{\mathfrak{N} \mathfrak{a}}& 2 d 
\frac{\Im r}{\mathfrak{N} \mathfrak{a}} \end{matrix}\right)\cdot
\left(\begin{matrix} \Re q & d\, \Im q \\ \Im q & \Re q\end{matrix}\right) 
\;\; &\text{if $d\not\equiv 1 \mo 4$},\\
\left(\begin{matrix} 2 & 0\\ 2\frac{\Re r}{\mathfrak{N} \mathfrak{a}}& 2 d 
\frac{\Im r}{\mathfrak{N} \mathfrak{a}} \end{matrix}\right)\cdot
\left(\begin{matrix} \Re q & d\, \Im q \\ \Im q & \Re q\end{matrix}\right) 
\cdot
\left(\begin{matrix} 1 & \frac{1}{2} \\ 0 & \frac{1}{2}\end{matrix}\right)
\;\; &\text{if $d \equiv 1 \mo 4$}.\end{aligned} \]

The matrix corresponding to $d\not\equiv 1$ has determinant
$4 d \,\Im r \mathfrak{N} q/\mathfrak{N} \mathfrak{a}$,
whereas the matrix corresponding to $d\equiv 1$ has determinant
$2 d \,\Im r \mathfrak{N} q/\mathfrak{N} \mathfrak{a}$.
Hence $\kappa(\jmath(L_{\mathfrak{a},r}))$ is
 either the empty set or a lattice $L'_{\mathfrak{a},r}$ 
of index
\[\lbrack \mathbb{Z}^2 : L'_{\mathfrak{a},r} \rbrack =
\begin{cases}
4 d \,\Im r \mathfrak{N} q  \lbrack \mathfrak{a} : L_{\mathfrak{a},r} \rbrack
&\text{if $d \not\equiv 1 \mo 4$,}\\
2 d \,\Im r \mathfrak{N} q  
\lbrack \mathfrak{a} : L_{\mathfrak{a},r} \rbrack
&\text{if $d \equiv 1 \mo 4$,}\end{cases}\]
and $\kappa(\jmath(S_{\mathfrak{a},r}))$ is a convex set $S'_{\mathfrak{a},r}$
contained in
\[\lbrack - 3 |d| \frac{\mathscr{S}(r)}{\mathfrak{N} \mathfrak{a}} \mathscr{S} (q) N',
 3 |d| \frac{\mathscr{S}(r)}{\mathfrak{N} \mathfrak{a}} \mathscr{S} (q) N' \rbrack^2,\]
which is contained in 
\[\lbrack - 3 |d| \mathscr{S}(q)
 \frac{n_{s 1} (\log N)^{\alpha+1}}{n_{s 0}},
3 |d| \mathscr{S}(q)
 \frac{n_{s 1} (\log N)^{\alpha+1}}{n_{s 0}}\rbrack^2 ,\]
which is in turn contained in
\[
\lbrack - k' (\log N)^{2 \alpha + 2} N, k' (\log N)^{2 \alpha + 2} N \rbrack^2,\]
where $k'$ depends only on $d$ and $q$. Write (\ref{eq:somename4}) as
\begin{equation}\label{eq:somename5}
\mathop{\mathop{\sum_{\mathfrak{a}}}_{n_{s 0} < \mathfrak{N}
 \mathfrak{a} \leq n_{s 1}}}_{n>1\Rightarrow n\nmid \mathfrak{a}}
\mathop{\sum_{r\in \overline{\mathfrak{a}} \cap R_s}}_{\Im r \ne 0}
\left|\sum_{(v,w)\in L'_{\mathfrak{a},r} \cap  S'_{\mathfrak{a},r}}
\lambda(v) \lambda(w)\right|
\end{equation}

Since $r$ is in $R_s$, $\Im r$ takes values between $-k n_{s 1} (\log N)^{\alpha+1}$
and $k n_{s 1} (\log N)^{\alpha+1}$. By Lemma \ref{lem:divbyno}, $\Im r$ takes each of
these
values  at most
\[2 \lceil (k n_{s 1} (\log N)^{\alpha+1})/n_{s 0} \rceil \ll
(\log N)^{2\alpha +2} \]
times. Thus (\ref{eq:somename5}) is bounded by a constant times
\[ \frac{N''}{2^{s-1}} (\log N)^{2\alpha + 2}
\sum_{0<y\leq k M_2 (\log N)^{\alpha +1}}
\max_\mathfrak{a} \max_{r : \Im r = y}
\left|\sum_{(v,w)\in L'_{\mathfrak{a},r} \cap  S'_{\mathfrak{a},y}}
\lambda(v) \lambda(w)\right|.\]
We may assume 
that $\lbrack \mathbb{Z}^2 : L\rbrack \ll (\log N)^\alpha$,
as otherwise what we want to prove is trivial. 
By Corollary \ref{bomb5},
\[\sum_{0<y\leq k M_2 (\log N)^{\alpha +1}}
\max_\mathfrak{a} \max_{r : \Im r = y}
\left|\sum_{(v,w)\in L'_{\mathfrak{a},r} \cap  S'_{\mathfrak{a},y}}
\lambda(v) \lambda(w)\right|\] is
\[O\left( \tau(4 d \num(\mathfrak{N} q) 
\den(\mathfrak{N} q) \lbrack\mathscr{O}_K : L'\rbrack^2)
\frac{((\log N)^{2\alpha +2} N)^2}{(\log N)^{8 \alpha + 9}}\right) ,\]
where $\den(a/b)$ stands for the denominator $b$
of a rational number $a/b$, $\gcd(a,b)=1$.
(Recall that $\num(a/b)$ stands for the numerator $a$.)
It is time to collect all terms:
\[\begin{aligned}
T_s^2 
&\ll 2^{s-1} \frac{N''}{2^{s-1}} (\log N)^{2\alpha+2}
\tau(4 d \num(\mathfrak{N} q) 
\det(\mathfrak{N} q) \lbrack \mathscr{O}_K : L'\rbrack^2)
\frac{((\log N)^{2\alpha +2} N)^2}{(\log N)^{8 \alpha+ 9}}\\
&+ 2^{s-1} \left(N^2 (\log N)^{2 (\alpha+1)} +
\frac{N^4 (\log N)^{3\alpha +3}}{2^s \sqrt{M_1}}\right)\\
&\ll N^2 (\log N)^{2 \alpha + 2 + \epsilon}
\frac{(\log N)^{4 \alpha + 4} N^2}{(\log N)^{8\alpha+9}} +
2^s N^2 (\log N)^{2 (\alpha+1)} +
N^4/(\log N)^{7\alpha+7}\\ &\ll N^4/(\log N)^{2\alpha +3 -\epsilon}
+ 2^s N^2 (\log N)^{2(\alpha+1)}
.\end{aligned}\]
Thus
\[\sum_{s=1}^{\lceil \log_2 (N''/M_1)\rceil} T_s \ll
\frac{N^2}{(\log N)^{\alpha+1}} + \frac{N^2 (\log N)^{\alpha+1}}{
(\log N)^{10 (\alpha+1)}} \ll
\frac{N^2}{(\log N)^{\alpha+1}} .\]
By (\ref{eq:satun})
and (\ref{eq:lomsuf}), the left-hand side of (\ref{eq:satun}) is
at most
\[O\left(\frac{\log M_1}{\log M_2} \frac{\Area(S^*)}{
\lbrack \mathscr{O}_K : L'\rbrack}+ N' M_2^2 \right)+
\sum_{s=1}^{\lceil \log_2 (N''/M_1)\rceil} T_s,\]
which, by the above, is at most
\[O\left(\frac{\log \log N}{\log N} \frac{\Area(S^*)}{
\lbrack \mathbb{Z}^2 : L\rbrack} + \frac{N^2}{(\log N)^{36 \alpha + 48}}
+ \frac{N^2}{(\log N)^{\alpha + 1}}\right)\]
As we saw before, there are at most $O(N^2/(\log N)^{\alpha})$
terms in the original sum missing from the left-hand side of
(\ref{eq:satun}). Hence, the original sum is at most
\[O\left(\frac{\log \log N}{\log N} \frac{\Area(S)}{
\lbrack \mathbb{Z}^2 : L\rbrack} +
\frac{N^2}{(\log N)^\alpha}\right),\]
as was to be proved.
\end{proof}
\section{Averages twisted by characters. Root numbers.}\label{sec:theman}
Define $(a|b) = \prod_{p|b, \text{$p$ odd}} (a/p)^{v_p(b)}$, where
$(a/p)$ is the quadratic reciprocity symbol.
\begin{prop}\label{prop:three}
Let $S$ be a convex subset of $\lbrack -N,N\rbrack^2$, $N>1$. Let $L\subset \mathbb{Z}^2$ be a lattice coset. Then
\[\mathop{\sum_{(x,y)\in S\cap L}}_{\gcd(x,y)=1} (y|x)\cdot 
\lambda(x y (x - y))
\ll \frac{\log \log N}{\sqrt{\log N}} \frac{\Area(S)}{
\lbrack \mathbb{Z}^2 : L\rbrack}
+ \frac{N^2}{(\log N)^\alpha}\]
for any $\alpha>0$. The implied constant depends only on $\alpha$.
\end{prop}
\begin{proof}
Proceed as in Theorem \ref{thm:one}. When applying Lemma \ref{split1},
use $M_1 = (\log N)^{2 \alpha + 2}$, $M_2 = e^{\epsilon \sqrt{\log N}}$,
where $\epsilon>0$ is a sufficiently small constant; define
$\mathscr{P}$ to be the set of all primes $p$ such that
$M_1\leq p<M_2$ and $p$ is not the largest prime divisor of
any exceptional modulus $m\leq N$.
Instead of (\ref{eq:coal1}),
we have
\[
\mathop{\sum_{a_1= s M_1}^{(s+1) M_1 -1}
\sum_{a_2= s M_1}^{(s+1) M_1 -1}}_{a_1\ne a_2} |\sigma_{a_1}
\sigma_{a_2}|
\left|\sum_{(b,c)\in S_{a_1,a_2}''\cap L_{a_1,a_2}''} 
(c|a_1 a_2) \lambda(a_1 b - c) \lambda(a_2 b - c) \right|.
\]
For every $k \mo a_1 a_2$, the lattice $L_{a_1,a_2,k}'''=
L_{a_1,a_2}'' \cap \{(b,c) : c\equiv k \mo a_1 a_2\}$ has index dividing
$\lbrack \mathbb{Z}^2 : L_{a_1,a_2}''\rbrack \cdot a_1 a_2$. Define
$L_{a_1,a_2,k} = \left(\begin{matrix}a_1 & -1\\a_2 &-1\end{matrix}\right)
L_{a_1,a_2,k}'''$, and proceed as before. We obtain, instead of 
(\ref{eq:coal2}),
\[M_1^2 \mathop{\max_{s M_1\leq a < (s+1) M_1}}_{\text{$a$
unexceptional}} 
    \mathop{\mathop{\max_{-M_1\leq d\leq M_1}}_{d\ne 0}}_{\text{$a+d$
unexceptional}}
 \sum_{k=1}^{a(a+d)} \left|
\sum_{(v,w)\in S_{a,a+d}\cap L_{a,a+d,k}} 
\lambda(v)
 \lambda(w) \right| .
\]
Since $a (a+d) \ll e^{2 \epsilon \sqrt{\log N}}$, the additional sum 
$\sum_k$ is absorbed by the right-hand side of the analogue of
(\ref{eq:coal3}). 
No invocation of Lemma \ref{lem:diverto} is needed here,
as the Siegel estimate (\ref{eq:sw}) suffices; in particular, no short-interval results
are needed. Note also that the factor of
$\prod_{M_1\leq p < M_2:\; p\notin \mathscr{P}}
\left( 1 - \frac{1}{p}\right)^{-1}$ coming from Lemma \ref{lem:sieve}
is bounded from above by a constant: since there are at most $O(\log \log N)$
exceptional moduli up to $N$, the cardinality of 
$\{M_1\leq p< M_2: p\notin \mathscr{P}\}$ is 
$O(\log \log N)$; as $M_1 = (\log N)^{2 \alpha + 2}$, it follows that
the factor
is $O(1)$. Note, finally, that
every time we apply the inequality (\ref{eq:sw}), the modulus $m$ is
not divisible by any exceptional moduli $q > \lbrack \mathbb{Z}^2 :
L_{a_1,a_2}''\rbrack$. Given that $\lbrack \mathbb{Z}^2 :
L_{a_1,a_2}''\rbrack \ll \lbrack \mathbb{Z}^2 : L\rbrack^2$ and that
we may assume that $\lbrack \mathbb{Z}^2 : L\rbrack < (\log N)^{\alpha}$
(as otherwise the result to be proven is trivial), we have that $m$
is not divisible by any exceptional moduli $q > C (\log N)^{2 \alpha}$,
where $C$ is an absolute constant, and thus (\ref{eq:sw}) is valid.
%The result follows -- with effectivity for $\alpha < 1$.
\end{proof}
\begin{cor}
For $a,b\in \mathbb{Z}$ coprime, let $E_{a,b}$ denote the elliptic curve
$y^2 = x (x + a) (x + b)$. Let $S$ be a convex subset of 
$\lbrack - N, N\rbrack^2$, $N>1$. Let $L\subset \mathbb{Z}^2$ be a lattice
coset. Then
\[\mathop{\sum_{(a,b)\in S \cap L}}_{\gcd(a,b)=1} 
 W(E_{a,b}) \ll \frac{\log \log N}{\sqrt{\log N}} \frac{\Area(S)}{\lbrack
\mathbb{Z}^2 : L\rbrack} + \frac{N^2}{(\log N)^{\alpha}}\]
for any $\alpha>0$, where $W(E_{a,b})$ is the root number of $E_{a,b}$. The implied constant depends only on $\alpha$.
\end{cor}
\begin{proof}
This is \cite{He}, Proposition 5.8. It is an easy consequence of Prop. 
\ref{prop:three}: the root number of $E_{a,b}$ equals
\[- \left(\frac{a}{\rad(2^{-v_2(b)} b)}\right)
\left(\frac{b}{\rad(2^{-v_2(a)} a)}\right)
\left(\frac{-a}{\rad(2^{-v_2(b-a)} (b-a))}\right)
\mu(\rad(a b (a - b))) ,\]
and the ratio of this expression to $(y|x) \lambda(x y(x-y))$ can be
handled by means of a square-free sieve (\cite{Hesq}, Prop. 3.12).
\end{proof}
\section{Acknowledgements}
The help and counsel of my Doktorvater, H. Iwaniec, were
 most valuable. Thanks are also due to an anonymous referee for his detailed
work.

\affiliationone{H. A. Helfgott\\
D\'epartement de Math\'ematiques et Statistique\\
Universit\'e de Montr\'eal\\CP 6128 succ Centre-Ville\\ 
Montr\'eal, QC H3C 3J7\\ Canada
\email{helfgott@dms.umontreal.ca}}
\end{document}